\documentclass[12pt]{amsart}
\usepackage{amsfonts}
\usepackage{latexsym}
\usepackage{amssymb}
\usepackage{amsmath}
\usepackage{amstext}
\usepackage{amsthm}
\usepackage{amsxtra}
\usepackage[utf8]{inputenc}
\usepackage[brazilian,english]{babel}
\providecommand{\U}[1]{\protect\rule{.1in}{.1in}}
\newtheorem{theorem}{Theorem}
\newtheorem{lemma}{Lemma}
\newtheorem{corollary}{Corollary}

\newtheorem{proposition}{Proposition}

\newtheorem{conjecture}{Conjecture}
\theoremstyle{definition}

\usepackage{enumerate}
\newtheorem*{theoremA}{Theorem A}

\begin{document}
	
	\title{On an internal characterization of horocyclically convex domains in the unit disk}
	
	\author[J. Arango]{Juan Arango}
	\address{Departamento de Matem\'aticas, Universidad Nacional de Colombia, Medell\'{\i}n, Colombia}
	\email{jharango@unal.edu.co}
	
	\author[H. Arbel\'aez]{Hugo Arbel\'{a}ez}
	\address{Departamento de Matem\'aticas, Universidad Nacional de Colombia, Medell\'{\i}n, Colombia}
	\email{hjarbela@unal.edu.co}
	
	\author[D. Mej\'ia]{Diego Mej\'ia}
	\address{Departamento de Matem\'aticas, Universidad Nacional de Colombia, Medell\'{\i}n, Colombia}
	\email{dmejia@unal.edu.co}

	\thanks{The first and third authors were supported by Universidad Nacional de Colombia, under the project Hermes, code 57923. The second author was supported by Universidad Nacional de Colombia, under the project Hermes, code 61126.}

	\maketitle

	\begin{abstract}
		A proper subdomain $G$ of the unit disk $\mathbb{D}$ is horocyclically convex (horo-convex) if, for every $\omega \in \mathbb{D}\cap \partial G$, there exists a horodisk $H$ such that $\omega \in \partial H$ and $G\cap H=\emptyset$. In this paper we give an internal characterization of these domains, namely, that $G$ is horo-convex if and only if any two points can be joined inside $G$ by a $C^1$ curve composed with finitely many Jordan arcs with hyperbolic curvature in $(-2,2)$. We also give a lower bound for the hyperbolic metric of horo-convex regions and some consequences.
	\end{abstract}
	
	
	
	\section{Introduction and main results}\label{section-1: introduction}

	A domain $G$ in the complex plane $\mathbb{C}$ is convex if and only if, for every $\omega \in \partial G$, there exists a half-plane $H$ such that $\omega \in \partial H$ and $G\cap H=\emptyset$. The line $\partial H$ is a noncompact maximal curve of constant euclidean curvature. This fact motivates to search for the analog within the hyperbolic geometry of the unit disk $\mathbb{D}$. The study of this analogy was initiated by Mej\'ia and Pommerenke in \cite{MePo4}.\\
	\indent Consider $\mathbb{D}$ endowed with the hyperbolic metric $\lambda_{\mathbb{D}}(z)|dz|=\vert dz \vert /(1-\vert z \vert^2)$. In this space the noncompact maximal curves of constant hyperbolic curvature $\kappa$ are those with $\vert \kappa \vert\leq 2$. These are circular arcs from $\mathbb{T}=\partial\mathbb{D}$ to $\mathbb{T}$. Hence, the analogy in the \emph{hyperbolic disk}, $(\mathbb{D},ds)$, of convexity in the euclidean plane can be established in multiple directions; among them, probably the most important is known as hyperbolic convexity: a domain $G\subset \mathbb{D}$ is hyperbolically convex (h-convex) if, for every $\omega \in \mathbb{D}\cap \partial G$, there is a hyperbolic half plane $H$ with $\omega \in \partial H$ and $G\cap H=\emptyset$. The supporting maximal curve at each point $\omega \in \mathbb{D}\cap \partial G$ has hyperbolic curvature zero. Hyperbolic convexity has been extensively studied; see for example [\cite{MaMi1}; \cite{MaMi2}; \cite{MePo1}; \cite{MePo2}; \cite{MePo3}; \cite{MePoVa}].\\
	\indent In \cite{MePo4} the authors considered the extremal case where the supporting maximal curve has hyperbolic curvature $\vert\kappa\vert =2$. A curve of this type touches $\mathbb{T}$ and is named \emph{horocycle}; the inner domain of a horocycle is a \emph{horodisk}.
	A domain $G\subset \mathbb{D}$ is called \emph{horocyclically convex} (horo-convex) if, for every $\omega \in \mathbb{D}\cap \partial G$, there exists a horodisk $H$ such that $\omega \in \partial H$ and $G\cap H=\emptyset$. A \emph{horocyclically convex function} $f$ is a conformal map of $\mathbb{D}$ onto a horo-convex domain $G\subset \mathbb{D}$.
	Clearly every h-convex domain is horo-convex, and every horo-convex domain is simply connected as we will show later.
	
	Let $\Gamma$ be a Jordan arc or a Jordan curve of the form $\Gamma = \Gamma_1 \cup \cdots \cup \Gamma_n$, where $\Gamma_k,\,k=1,\dots,n$, are smooth Jordan  arcs from $p_{k-1}$ to $p_k$, that are otherwise disjoint. We denote by $\Delta(\Gamma_k)$ the change of the tangent angle along $\Gamma_k$, and by $\Delta(p_k),\,k=1,\dots,n-1$, the change of the tangent angle at the vertex $p_k$. Then, the change of the tangent angle along $\Gamma$ is given by
\begin{equation}\label{eq-1.1}
	\Delta(\Gamma)=\sum_{k=1}^n\Delta(\Gamma_k)+\sum_{k=1}^{n-1}\Delta(p_k).
\end{equation}
If $\Gamma$ is a positively oriented smooth Jordan curve, then $\Delta(\Gamma)=2\pi$. We are concerned with the case where the arcs $\Gamma_k$ are all inside $\mathbb{D}$, have constant hyperbolic curvature in the open interval $(-2,2)$, and in such a way that $\Gamma_{k-1}$ and $\Gamma_k$ join smoothly at the point in common, making $\Gamma$ of class $C^1$. We call \emph{admissible} a Jordan arc of this kind.

Our main objective is to find internal characterizations of horo-convexity. In this direction we have the following result.
\vspace{.2cm}
\begin{theorem}\label{internalCh}
	Let $G\subset\mathbb{D}$ be a simply connected domain with locally connected boundary. The following conditions are equivalent,
	\begin{enumerate}
		\item[(i)] $G$ is horocyclically convex.
		\item[(ii)] For every $a,b\in G$, with $a\neq b$, there is an admissible Jordan arc $\Gamma \subset G$ joining $a$ and $b$.
	\end{enumerate} 
\end{theorem}

The condition (ii) on the Jordan arc $\Gamma$ does not imply simple connectivity. Take for instance the unit disk punctured at the origin, $\mathbb{D}^*$: any two points in $\mathbb{D}^*$ can be joined inside $\mathbb{D}^*$ by an arc with constant hyperbolic curvature between $0$ and $2$. On the other hand, it is possible that the imposed condition of local connectivity on the boundary, be superfluous.\\\\
\indent Another topic in which we are interested deal with finding a characterization of horo-convexity by means of a lower bound for the hyperbolic metric in terms of the hyperbolic distance to the boundary. This type of characterization is known for convex domains in all three geometries: euclidean (\cite{MeMi1}), hyperbolic (\cite{Su2}) and spherical (\cite{Su1}). We recall that the hyperbolic distance in the unit disk is given by
\[
d_{\mathbb{D}}(z,w)=\operatorname{arctanh}\bigg \vert \frac{z-w}{1-\bar{z}w}\bigg \vert=\operatorname{arctanh}e_{\mathbb{D}}(z,w),
\]
where
\[
e_{\mathbb{D}}(z,w)=\bigg \vert \frac{z-w}{1-\bar{z}w}\bigg \vert
\]
is the so called pseudo-hyperbolic distance in $\mathbb{D}$. Since $\operatorname{arctanh}x\geq x$ for $0\leq x <1$, we have the inequality $d_{\mathbb{D}}(z,w)\geq e_{\mathbb{D}}(z,w)$.\\ For a proper subdomain $G$ of $\mathbb{D}$ we denote by $\lambda_G(z)|dz|$ its hyperbolic metric and by $\nu_G(z)$ its hyperbolic density; that is,
\[
\nu_G(z)=\frac{\lambda_G(z)|dz|}{\lambda_{\mathbb{D}}(z)|dz|}=(1-|z|^2)\lambda_G(z).
\]

The hyperbolic density is a continuous function on $G$ which is invariant under the group \textrm{M\"ob}($\mathbb{D}$) of conformal automorphisms of $\mathbb{D}$. Also, for $z\in G$, we denote by $d_G(z)$ and $e_G(z)$, respectively, the hyperbolic and pseudo-hyperbolic distance from $z$ to the boundary of $G$. As in \cite{Su2}, we consider the following quantities
\[
C_{\mathbb{D}}(G):=\inf_{z\in G}d_G(z)\nu_G(z),\qquad C'_{\mathbb{D}}(G):=\inf_{z\in G}e_G(z)\nu_G(z).
\]
Clearly, $C_{\mathbb{D}}(G)\geq C'_{\mathbb{D}}(G)$.\\\\
 We name by \emph{horo-crescent} the standard horo-convex domain which is the exterior in $\mathbb{D}$ of a closed horodisk. Since any two horodisks are \textrm{M\"ob}($\mathbb{D}$)-equivalent, then any two horo-crescent domains are also \textrm{M\"ob}($\mathbb{D}$)-equivalent.\\
 At this point we refer the reader to \cite{BeMi} for fundamental facts of the hyperbolic metric.  Based on the definition of horo-convexity and the monotonicity property of the hyperbolic metric we have,
\begin{theorem}\label{lower_bound}
	Let $G\subset\mathbb{D}$ be any horo-convex domain. Then, for each $z\in G$,
	\begin{equation}\label{eq-1.2}
		\nu_G(z)\geq\frac{\pi}{e^{2d_G(z)}}\frac{1}{\sin \dfrac{\pi}{e^{2d_G(z)}}},
	\end{equation}
	with equality at a point if and only if $G$ is a horo-crescent domain.
\end{theorem}

A natural question to ask is for what kind of simply connected regions in the unit disk is valid inequality \eqref{eq-1.2}. In this direction one can also ask if it might be possible to extend the method developed by Sugawa in (\cite{Su1}, \cite{Su2}) to horocyclical convexity. We don't have answers to these questions. Sugawa gave versions of Keogh's lemma (\cite{Ke}) to spherical and hyperbolic convexity which allowed to demonstrate that lower bounds of the hyperbolic metric for spherical and hyperbolic convex domains indeed characterize those types of domains. \\\\
For $t>0$ define $h$ by
\[
h(t):= \frac{e^{2t}}{\pi}\sin\frac{\pi}{e^{2t}}.
\]
Using the software Geogebra we see that the minimum value of $g(t)=t/h(t)$ is $\approx 0.48$ and occurs as $\approx 0.11$; so $C_{\mathbb{D}}(G)\gtrapprox 0.48$.\\\
On the other hand, in terms of $e_G(z)$ inequality \eqref{eq-1.2} is 
\begin{equation}\label{eq-1.3}
	\nu_G(z)\geq \frac{1-e_G(z)}{1+e_G(z)}\frac{\pi}{\sin \dfrac{\pi(1-e_G(z))}{1+e_G(z)}};
\end{equation}  
 again, using Geogebra we see that the minimum value of
\[
\frac{1-s}{1+s}\frac{\pi s}{\sin \dfrac{\pi(1-s)}{1+s}}, \quad 0<s<1,
\]
is $\approx 0.48$ and occurs at $\approx 0.12$. So $C'_{\mathbb{D}}(G)\gtrapprox 0.48$. In view of the results proved by Sugawa mentioned above, it is reasonable to think that any simply connected domain $G$ in $\mathbb{D}$ for which $C'_{\mathbb{D}}(G)< 0.48$ is not horo-convex. With this in mind we define the class $\mathcal{C}$ of domains $\Omega\subset \mathbb{D}$ of the form $D_1\setminus\overline{D_2}$ where each $D_j$ is a hyperbolic disk in $\mathbb{D}$ and such that $\partial D_1$ and $\partial D_2$ intersect orthogonally at two points in $\mathbb{D}$. Let $m(\Omega)$ denote the hyperbolic midpoint of the hyperbollically concave boundary arc $\partial D_2\cap D_1$ of $\Omega$. We recall that the hyperbolic curvature of $\partial D_j,\;j=1,2$, is bigger than two.\\
\begin{conjecture}
	Let $G$ be a simply connected subdomain of $\mathbb{D}$. Then, $G$ is not horociclically convex in $\mathbb{D}$ if and only if there is a domain $\Omega \in \mathcal{C}$ such that $\Omega \subset G$ with $m(\Omega)\in\partial G\cap \mathbb{D}$.
\end{conjecture} 
The Principle of hyperbolic metric establishes that if $G$ is a simply connected subregion of $\mathbb{D}$ and $f$ is holomorphic on $\mathbb{D}$ with $f(\mathbb{D})\subset G$, then $$\nu_G(f(z))\frac{|f'(z)|}{1-|f(z)|^2}\leq \lambda_{\mathbb{D}}(z)$$
for $z\in \mathbb{D}$ with equality if and only if $f$ is a conformal mapping of $\mathbb{D}$ onto $G$. 
As a consequence of Theorem \ref{internalCh} and the Principle of hyperbolic metric is the following corollary.
\begin{corollary}\label{corollary_lowerbound}
	Let $G\subset\mathbb{D}$ be a horo-convex domain and $f:\mathbb{D}\to \mathbb{D}$ a holomorphic mapping with $f(\mathbb{D})\subset G$. Then, for $z\in \mathbb{D}$,
	\begin{equation}\label{eq-1.4}
		(1-|z|^2)\frac{|f'(z)|}{1-|f(z)|^2}\leq h(d_G(f(z))).
	\end{equation}
Equality holds at a point if and only if $G$ is a horo-crescent domain and $f$ is a conformal mapping of $\mathbb{D}$ onto $G$. 
\end{corollary}
Since $h$ is strictly increasing, by applying \eqref{eq-1.4} at $z=0$ we see that the range of $f$ covers de hyperbolic disk centered at $f(0)$ with hyperbolic radius $h^{-1}(|f'(0)|/(1-|f(0)|^2))$.

	\section{Proofs of results }\label{section-2:proofs}
	\subsection{Simple connectivity of horo-convex domains}\label{sec1:subsec1} We need the following lemma.
	\begin{lemma}\label{sec2:subsec2.1:lemma1}
		Let $G$ be a domain contained in $\mathbb{D}$. If $B$ is a component of $\hat{\mathbb{C}}\setminus G$ with $\infty \notin B$, then $\partial B\subset \partial G$.
	\end{lemma}
	
	\begin{proof}
		Let $B$ be a component of $\hat{\mathbb{C}}\setminus G$ as in the statement, and $x\in \partial B$. Since $B$ is closed in $\hat{\mathbb{C}}$ then $x\notin G$. Since $G$ is open in $\hat{\mathbb{C}}$ then $\hat{\mathbb{C}}=G\cup \partial G\cup \operatorname{Ext}(G)$. Suppose that $\infty \neq x\in \operatorname{Ext}(G)$. Then, there exists an euclidean open disk $D$ centered at $x$ such that $D\subset \operatorname{Ext}(G)$. Hence $D\subset B$. So, $x\in \operatorname{Int}(B)$. This is a contradiction. Therefore $x\in \partial G$.
	\end{proof}
	\begin{proposition}\label{sec:subsec2.1:lemma1}
		Every horo-convex domain is simply connected.
	\end{proposition}
	\begin{proof}
		Let $G\subset \mathbb{D}$ be horo-convex. We must prove that $\hat{\mathbb{C}}\setminus G$ is connected. Note that $\hat{\mathbb{C}}\setminus \mathbb{D}\subset\hat{\mathbb{C}}\setminus G$ and $\hat{\mathbb{C}}\setminus \mathbb{D}$ is connected. Let $A$ be the component of $\hat{\mathbb{C}}\setminus G$ that contains $\hat{\mathbb{C}}\setminus \mathbb{D}$. We argue by contradiction assuming that $\hat{\mathbb{C}}\setminus G$ is not connected. Then, there exists a component $B$ of $\hat{\mathbb{C}}\setminus G$ different from $A$. Hence $B\subset \hat{\mathbb{C}}\setminus A\subset \mathbb{D} $. Since $B$ is closed in $\hat{\mathbb{C}}\setminus G$ which is closed in $\hat{\mathbb{C}}$, then $B$ is closed in $\hat{\mathbb{C}}$. Furthermore, the boundary of $B$, $\partial B$, respect to $\hat{\mathbb{C}}$ is not empty (otherwise, $B=\overline{B}=B^\circ\cup \partial B = B^\circ$; then $B$ would be open and closed in $\hat{\mathbb{C}}$, contradicting the fact that $\hat{\mathbb{C}}\setminus G$ is a proper subset of $\hat{\mathbb{C}}$). Let $\zeta \in \partial B$. So, $\zeta \in \overline{B}=B\subset\mathbb{D}$. By Lemma \ref{sec2:subsec2.1:lemma1}, $\zeta \in \partial G \cap \mathbb{D}$. By the horo-convexity of $G$ there exists a horo-disk $H$ with $\zeta \in \partial H \cap \mathbb{D}$ and $G\subset \mathbb{D}\setminus \overline{H}$. Since $\zeta \in B\cap \overline{H}$, then $B\cup \overline{H}$ is connected in $\hat{\mathbb{C}}\setminus G$. Therefore $B\cup \overline{H}$ lies inside a unique component of $\hat{\mathbb{C}}\setminus G$; hence $B\cup \overline{H}\subset B$ which implies $\overline{H}\subset B$. This gives a contradiction since $\overline{H}\cap \mathbb{T}\neq \emptyset$ and $\mathbb{T}\subset A$.
	\end{proof}
	\subsection{Internal characterization of horo-convex domains}\label{sec1:subsec2}
	
	The proof of Theorem \ref{internalCh} is very geometric and depends on several lemmas some of which are inspired on ideas already present in \cite{MePo4}. Also, sometimes we prefer to use the upper half-plane $\mathbb{H}=\{z\in\mathbb{C}: \operatorname{Im}z>0\}$ as our model for the hyperbolic space. We recall that in this model the arcs of constant hyperbolic curvature in $(-2,2)$ meet the boundary of $\mathbb{H}$.
	\begin{lemma}\label{sec2:subsec2.2:lemma2}
		Let $H\subset \mathbb{D}$ be any horodisk.
		\begin{enumerate}
			\item[(i)] If $\Gamma$ is an admissible Jordan arc leaving $H$ at $w\in\partial H \cap \mathbb{D}$, then $\Gamma$ cannot come back to $\partial H$.
			\item[(ii)] If $a,b$ are two different points in $H$ and $\Gamma$ is an admissible Jordan arc with end points $a,b$, then $\Gamma$ is contained in $H$. 
		\end{enumerate}  
	\end{lemma}
	\begin{proof}
		Since the hyperbolic curvature is invariant under conformal mappings, we may work in the upper half-plane $\mathbb{H}$, and we may also suppose that $H=\{z: \operatorname{Im}z>c\}$ for some $c>0$.\\
		$(i)$ By hypothesis, the initial direction of $\Gamma$ at $w$ points downward. If $\Gamma $ turns back to $\partial H$, it would exist a point $p\in \Gamma$ with minimal imaginary part. Let $\Gamma_k$ be an arc component of $\Gamma$ containing $p$. Since $\Gamma$ is of class $C^1$, the tangent line to $\Gamma_k$ at $p$ has to be horizontal. But any arc in $\mathbb{H}$ with constant hyperbolic curvature in $(-2,2)$ can only have horizontal tangent at a point with maximal imaginary part. This is a contradiction.\\
		$(ii)$ By part $(i)$, $\Gamma$ cannot leaves out $H$ at any point of $\partial H$, so $\Gamma$ stays in $\overline{H}$. Now, if $\Gamma$ meets $\partial H$ at a point $w$, then there would be an arc $B$, component of $\Gamma$, with $w\in B$. Since $\Gamma\subset \overline{H}$ is of class $C^1$, then $\partial H$ is tangent to $B$ at $w$. Therefore $\operatorname{Im}z,\,z\in B$, will have an strict minimum at $w$. But this is not possible because the hyperbolic curvature of $B$ belongs to $(-2,2)$. 
	\end{proof}
	
	Let $a,b$ be different points in $\mathbb{D}$. Let $H^+$, $H^-$ be the unique (open) horodisks whose boundaries pass through $a,b$. We denote $E(a,b):=H^+\cap H^-$.
	\begin{lemma}\label{sec2:subsec2.2:lemma3}
		Let $a,b$ be different points in $\mathbb{D}$. If $\Gamma$ is an admissible Jordan arc with end points $a$ and $b$, then $\Gamma \subset E(a,b)$.
	\end{lemma}
	\begin{proof}
		This is an immediate consequence of the previous lemma.
	\end{proof}	
	\begin{lemma}\label{sec2:subsec2.2:lemma4}
		Let $G\subset\mathbb{D}$ be a simply connected domain. If any two points of $G$ can be joined inside $G$ by an admissible Jordan arc $\Gamma$, then $\vert \Delta(\Gamma) \vert \leq 8\pi$. 
	\end{lemma}
	\begin{proof}
		Let $w_0,w_1 \in G$. We join $w_0$ and $w_1$ by an arc $\Gamma$ that satisfies the condition of the statement. Let $H_j^{\pm}$, $j=0,1$, be the horodisks tangent to $\Gamma$ at $w_j$ and touches $\mathbb{T}$ at $\zeta_j^{\pm}$. Since the horocycles have h-curvature with absolute value $2$, then $H_j^{\pm}\cap \Gamma=\emptyset$. We now construct a smooth positively oriented Jordan curve $J$ as follows. We may assume that $\Gamma$ is positively oriented from $w_0$ to $w_1$ and that $H_j^+$ lies above $\Gamma$ at $w_j$. Then $J=J_0 \cup \Gamma \cup J_1 \cup J^*$, where $J_j$ is the portion of $\partial H_j^+$ from $w_j$ to $\zeta_j^+$ that makes a $\pi$ angle with $\Gamma$ at $w_j$, and $J^*$ is the portion of $\mathbb{T}$ from $\zeta_1^+$ to $\zeta_0^+$. So,
		\begin{equation}\label{eq-2.2-1}
			2\pi = \Delta(J_0)+\Delta (\Gamma)+\Delta (J_1)+\Delta (J^*).
		\end{equation}
		It follows that
		\begin{equation}\label{eq-2.2-2}
			\vert\Delta(\Gamma)\vert \leq 8\pi.
		\end{equation}
	\end{proof}
	
	\begin{lemma}\label{sec2:subsec2.2:lemma5}
		If $f$ is any conformal map from $\mathbb{D}$ onto a domain $G$ that satisfies the condition (ii) of Theorem \ref{internalCh}, then $f'$ is in the Hardy class $H^p$ for some $p>0$. In particular, the angular derivative of $f$ exists and is different from $0,\infty$ almost everywhere.  
	\end{lemma}
	\begin{proof}
		Fix $w_0\in G$ and let $w$ be any other point in $G$. Consider a conformal map $f$ from $\mathbb{D}$ onto $G$ with $f(0)=w_0$ and $f(z)=w$. Let $C=f([0,z])$ and $\theta_j$, $j=0,1$, the angles between $[0,z]$ and $f^{-1}(\Gamma)$ at $0$ and $z$, respectively. Then,
		\begin{equation}\label{eq-2.2-3}
			\arg f'(z)-\arg f'(0)= \Delta (\Gamma) \pm \theta_0 \pm \theta_1.
		\end{equation}
		Therefore, it follows from (\ref{eq-2.2-2}) and (\ref{eq-2.2-3}) that
		\begin{equation}\label{eq-2.2-4}
			\vert \arg f'(z)-\arg f'(0)\vert \leq 12\pi.
		\end{equation}
		Inequality (\ref{eq-2.2-4}) implies that $f'$ is subordinate to the function
		\[
		g(z)=f'(0)\exp{\bigg[14\log\frac{1+z}{1-z}\bigg ] }=f'(0)\bigg (\frac{1+z}{1-z}\bigg )^{14}.
		\]
		Since $g\in H^p$ for $0<p<1/14$ the result follows from Littlewood's Subordination Principle (see \cite{Book_Duren}, p. 10).
	\end{proof}
	We will use the following geometric characterization of isogonality (see \cite{Book_Pommerenke}, p.254).
	\begin{theoremA}\label{isogonality}
		A conformal map $f$ from $\mathbb{D}$ onto a domain $G\subset \mathbb{C}$ is isogonal at $\zeta \in \mathbb{T}$ if and only if there is a curve $C\subset \mathbb{D}$ ending at $\zeta$ such that
		\[
		f(C):\;\omega+te^{i\alpha}, \quad 0<t\leq t_1
		\]
		satisfies
		\[
		\{w:\operatorname{Re}\big [e^{-i\alpha}(w-\omega)\big ]>\varepsilon t, \vert w-\omega\vert <t/\varepsilon\}\subset G
		\]
		for $0<\varepsilon<1$, $0<t<t_0(\varepsilon)$; and
		\[
		\text{there are points}\;\omega_t^{\pm}\in\partial G\;\text{with}\;e^{-i\alpha}(\omega_t^{\pm}-\omega)\sim \pm it
		\]
		as $t\to 0$.
	\end{theoremA} 
	Let $f$ be as in the previous lemma and let $$E=\{\zeta\in \mathbb{T}:\text{the angular derivative}\; f'(\zeta)\,\text{exists and is different from}\; 0,\infty\}.$$
	\begin{proof}[\textbf{Proof of $(ii)\implies (i)$ in Theorem \ref{internalCh}}]
		Since the boundary of $G$ is locally connected, $f$ can be extended continuously to $\overline{\mathbb{D}}$.\\
		$(A)$ We will show first that if $\zeta \in E$ with $f(\zeta)\in\mathbb{D}$, then there is a horodisk $H$ such that $f(\zeta)\in\partial H$ and $H\cap G=\emptyset$. 
		We may take $f(\zeta)=0$. By Theorem A there exists a curve $C\subset \mathbb{D}$ ending at $\zeta$ such that $f(C)=\{te^{i\alpha}:0<t\leq t_1\}$ for some $t_1>0$ and some $\alpha\in [0,2\pi)$. Furthermore, there exist $w_t^{\pm }\in\partial G$ such that 
		\[
		\lim_{t\to 0^+}\frac{e^{-i\alpha}w_t^{\pm}}{\pm it}=1 \;\text{and}\;\lim_{t\to 0^+}\vert w_t^{\pm}\vert=0.
		\]
		By rotating about the origin if necessary, we may take $\alpha=\pi/2$. So, $\lim_{t\to 0^+}\arg w_t^+=0$, and $\lim_{t\to 0^+}w_t^-=\pi$.\\
		Let $H$ be the horodisk such that $0\in\partial H$ and $\partial H$ is orthogonal to $f(C)$ at the origin. We will show that $H\cap G=\emptyset$. Arguing by contradiction suppose there is $q\in H\cap G$. Take a sequence $\{v_n\}$ in $f(C)$ converging to $0$ with $\vert v_{n+1}\vert < \vert v_n\vert$. By hypothesis, for each $n$ there is an admissible Jordan arc $\Gamma_n\subset G$ joining $v_n$ and $q$ . Let $\partial H^+=\partial H\cap\{z:\operatorname{Re}z>0 \}$ and $\partial H^-=\partial H\cap\{z:\operatorname{Re}z<0 \}$. We have two cases:\\
		$(a)$ There are $m,n$ such that $\Gamma_m$ gets in $H$ through $\partial H^+$ and $\Gamma_n$ gets in $H$ through $\partial H^-$. By Lemma \ref{sec2:subsec2.2:lemma2} there are arcs of $\Gamma_m$, $\Gamma_n$ and (possibly) of $f(C)$ that enclose a region $R$ with $0\in R$. Since $\Gamma_m$, $\Gamma_n$ and $f(C)$ are inside $G$, and $G$ is simply connected, then $R\subset G$, and so $0\in G$. This contradicts the fact that $0\in\partial G$.\\
		$(b)$ For all $n$, $\Gamma_n$ gets in $H$ through $\partial H^+$ (analogously, if for all $n$,  $\Gamma_n$ gets in $H$ through $\partial H^-$). There are two subcases:\\
		$(i)$ First consider $\operatorname{Re}q>0$. Let $O_0$ the horodisk such that $0,q\in\partial O_0$ and $\mathbb{R}^+\cap O_0=\emptyset$, and let $q_0$ be the first point (from $v_1$) where $\Gamma_1$ meets $\partial O_0$ (this point could be $q$ itself). Let $A$ be the arc of $\partial O_0$ between $0$ and $q_0$ open at $0$. We will show that $A\subset G$. Let $z\in A$. Take an open disk $D\subset G$ centered at $q_0$; then there exist $q_1\in D$ and $w\in f(C)$ with $\vert w\vert<\vert v_1 \vert$, such that the horodisk $O_1$, whose boundary passes through $q_1$ and $z$, also contains $w$. By hypothesis, there is an admissible Jordan arc $\Gamma \subset G$ joining $w$ and $q_1$. By Lemma \ref{sec2:subsec2.2:lemma3}, $\Gamma \subset E(w,q_1)\subset O_1$. If $B$ is the arc of $\Gamma_1$ between $q_0$ and $v_1$, then the closed curve $\Gamma \cup \overline{q_1q_0} \cup B \cup \overline{v_1w}$ enclose a region $R_1$. Since $G$ is simply connected, $R_1\subset G$. But, clearly, $z\in R_1$; so $z\in G$.\\
		Now let $\{te^{i\beta}:t\in\mathbb{R}\}$ be the tangent line to $\partial O_0$ at $0$; then $\beta \in (-\pi/2,0)$. There is $t>0$ such that $\vert w_t^+\vert < \operatorname{dist}(0,\Gamma_1)$ and $\vert \arg w_t^+\vert <\vert \beta \vert$. So, there exist $w\in f(C)$ and $z\in A$ such that if $O_2$ is a horodisk $w,z\in \partial O_2$, then $w_t^+\notin O_2$. By hypothesis there exists an admissible Jordan arc $\Gamma_0\subset G$ joining $w$ and $z$. Then, according to Lemma \ref{sec2:subsec2.2:lemma3}, $\Gamma_0\subset E(w,z)\subset O_2$. Finally, if $A'$ is the arc in $A$ between $z$ and $q_0$, then $J=\Gamma_0 \cup A' \cup B \cup \overline{v_1w}$ is a Jordan curve inside $G$. Since $G$ is simply connected, the interior of $J$, $R_0$, is contained in $G$. But $w_t^+\in R_0$, so $w_t^+\in G\cap \partial G$ which is not possible.\\
		$(ii)$ Now consider $\operatorname{Re}q \leq 0$. By part $(i)$, $H\cap G \cap \{\operatorname{Re}z>0\}=\emptyset$; hence, it is no possible for the $\Gamma_n$ to get in $H$ through $\partial H^+$.\\
		$(B)$ Now take $w_0\in \partial G\cap \mathbb{D}$. Then, there is $\zeta_0\in \mathbb{T}$ such that $f(\zeta_0)=w_0$. Since $\mathbb{T}\setminus E$ has measure $0$, $E$ is dense in $\mathbb{T}$. Therefore there exists a sequence $\{\zeta_n\}$ in $E$ such that $\zeta_0=\lim_{n\to\infty}\zeta_n$. By part $(A)$, for each $n$ there is a horodisk $H_n$ with $f(\zeta_n) \in \partial G$ and $H_n \cap G=\emptyset$. Let $c_n$ be the euclidean center of $H_n$. We may suppose that $\{c_n\}$ converges to $c_0\in\overline{\mathbb{D}}$. Indeed $c_0\in \mathbb{D}$ otherwise,  $r_n=|f(\zeta_n)-c_n|\to |f(\zeta_0)-c_0|=0$, contradicting the assumption $f(\zeta_0)=w_0\in\mathbb{D}$. So, we may suppose that $\{r_n\}$ converges to a positive number $r_0$. Let $H_0$ be the open disk centered at $c_0$ and radius $r_0$. Since $\operatorname{dist}(c_0,\mathbb{T})=\lim_{n\to\infty}\operatorname{dist}(c_n,\mathbb{T})=r_0$, then $H_0$ is a horodisk, and moreover, $w_0\in\partial H_0$. It only remains to show that $H_0\cap G=\emptyset$, but this is clear since $H_n\cap G=\emptyset$ for all $n$ and $H_n$ approaches $H_0$ as $n$ goes to $\infty$.
	\end{proof}
\begin{proof}[\textbf{Proof of $(i)\implies (ii)$ in Theorem \ref{internalCh}}] 
	Let $a,b\in G$ be given, with $a\neq b$. Take a Jordan arc $\Gamma_0\subset G$ joining $a$ and $b$. We cover $\partial G$ with open disks as follows: for each $\omega\in \partial G\cap \mathbb{D}$ there is, by hypothesis, a horodisk $H_{\omega}$ with $\omega \in \partial H_{\omega}$ and $H_{\omega}\cap G=\emptyset$. Now, take a slightly bigger euclidean open disk $D_{\omega}$ away from $\Gamma_0$, with $\omega\in D_{\omega} \supset H_{\omega}$; on the other hand, for each $\omega \in \partial G\cap\mathbb{T}$, take a small open euclidean disk $D_{\omega}$ away from $\Gamma_0$ , centered at $\omega$. By compactness, cover $\partial G$ with a finite number of such disks, $D_1,\dots,D_{\mu}$. Observe that, for $\omega\in\partial G$,  $\partial D_{\omega}\cap\mathbb{D}$ has hyperbolic curvature in $(-2,2)$.
	Changing to the upper half-plane model $\mathbb{H}$ where the argument is easier to visualize, we may suppose that the disks $D_i$ have euclidean center $c_i$ and euclidean radius $R_i$, with $0\leq \operatorname{Im}c_i<R_i$, $\partial G\subset \cup_{i=1}^{\mu}D_i$ and $\Gamma_0\cap \cup_{i=1}^{\mu}\overline{D_i}=\emptyset$. Also, we order the $D_i$'s in such a way that $\operatorname{Re}c_i\leq \operatorname{Re}c_{i+1}$. By invariance of the hyperbolic curvature, $\partial D_i\cap \mathbb{H}$ has hyperbolic curvature in $(-2,2)$ for each $i$. Notice also that non horizontal euclidean line segments have hyperbolic curvature in $(-2,2$). If the euclidean line segment $\overline{ab}$ is contained in $G$, take $\Gamma=\overline{ab}$ in case is not horizontal, otherwise, let $\Gamma$ be the arc between $a$ and $b$ of any circle  through $a$ and $b$ contained in $G$ that meets $\partial \mathbb{H}$.\\
	If $\overline{ab}$ is not contained in $G$ we have two general cases: $\overline{ab}$ is not vertical and $\overline{ab}$ is vertical.\\
	\textbf{A}: $\overline{ab}$ is not vertical. We may suppose without loss of generality that $\operatorname{Re}a<\operatorname{Re}b$. Let $J=\{1,2,\dots,\mu\}$ and $J_1=\{i\in J:\operatorname{Re}a\leq \operatorname{Re}c_i\leq \operatorname{Re}b, \overline{ab}\cap D_i\neq\emptyset\}$. We have two cases: $J_1\neq \emptyset$ and $J_1 = \emptyset$.\\
	\underline{$J_1\neq \emptyset$}: Let $j_1:=\min J_1$ and $h_1:=\max J_1$. For each $i\in J_1$ let $r_i\in \partial D_i$ be such that $\overline{ar_i}$ is the upper tangent segment from $a$ to $\partial D_i$, and let $\alpha_i:=\arg (r_i-a) \in (-\pi/2,\pi/2)$, $\alpha =\max\{\alpha_i:i\in J_1\}$, and $k:=\max\{i\in J_1:\alpha_i=\alpha\}$. Let $s\in \partial D_k$ be such that $\overline{bs}$ is the upper tangent segment from $b$ to $\partial D_k$. If $J_2:=\{i\in J_1:i>k \;\text{and}\;\overline{bs}\cap D_i\neq \emptyset\}=\emptyset$, take $\Gamma_1:=\overline{ar_k}\cup A\cup \overline{sb}$, where $A$ is the arc on $\partial D_k$ from $r_k$ to $s$. If $J_2\neq \emptyset$ let, for each $i\in J_2$, $s_i\in \partial D_k$ and $t_i\in\partial D_i$ be such that $\overline{s_it_i}$ is the upper tangent segment common to $\partial D_k$ and $\partial D_i$, and let $\beta_i:\arg (t_i-s_i)$ (which is less than $\alpha$). Let $\beta :=\max\{\beta_i: i\in J_2\}$ and $l:=\max\{i\in J_2:\beta_i=\beta\}$. Let $t\in \partial D_l$ be such that $\overline{bt}$ is the upper tangent segment from $b$ to $\partial D_l$. If $J_3:=\{i\in J_2:i>l \;\text{and}\;\overline{bt}\cap D_i\neq \emptyset\}=\emptyset$, take $\Gamma_1:=\overline{ar_k}\cup A'\cup \overline{s_lt_l}\cup B\cup \overline{tb}$, where $A'$is the arc on $\partial D_k$ from $r_k$ to $s_l$ and $B$ is the arc on $\partial D_l$ from $t_l$ to $t$. If $J_3\neq \emptyset$ we repeat the process. This process must end since $J$ is finite, and so we have constructed a ``polygonal" curve $\Gamma_1$ (with decreasing slope) from $a$ to $b$ composed by arcs of hyperbolic curvature in the interval $(-2,2)$. If $\Gamma_1\subset G$ take $\Gamma:=\Gamma_1$.\\
	Now we must take into consideration the possibility that $\Gamma_1\nsubseteq G$. Since the circular arcs of $\Gamma_1$ are indeed inside $G$, then $\Gamma_1$ fails to be inside $G$ when some (or all) its line segments meet $\partial G$ and, therefore, some $D_i$. By our construction, these indices $i$ must be either less than $j_1$ or bigger than $h_1$ and, for each one of theme the curve $\Gamma_1$, that we called of ``first level", has to be corrected by constructing ``polygonal" curves of ``second level" for each tangent segment that meets $\partial G$. We illustrate this as follows:\\
	Suppose for instance that $\overline{s_lt_l}\nsubseteq G$ and $\arg(t_l-s_l)>0$. Then $K_1:=\{j\in J:j<j_1\;\text{and}\;\overline{s_lt_l}\cap D_j\neq \emptyset\}\neq \emptyset$. For each $j\in K_1$ let $u_j\in \partial D_k$ and $v_j\in \partial D_j$ be such that $\overline{u_jv_j}$ is the upper cross tangent segment common to $\partial D_k$ and $\partial D_j$, and let $\gamma_j:=\arg (v_j-u_j)$. Let $\gamma:=\min \{\gamma_j:j\in K_1\}$ and $m:=\min \{j\in K_1:\gamma_j=\gamma\}$. Let $v\in \partial D_m$ and $t\in \partial D_l$ be such that $\overline{vt}$ is the cross tangent segment common to $\partial D_m$ and $\partial D_l$. If $K_2:=\{j\in K_1: j < m\;\text{and}\; \overline{vt}\cap D_j\neq \emptyset\}=\emptyset$, take $\Gamma_2:= C\cup \overline{u_mv_m}\cup D\cup \overline{vt} \cup E$, where $C$ is the arc on $\partial D_k$ from $s_i$ to $u_m$, $D$ is the arc on $\partial D_m$ from $v_m$ to $v$ and $E$ is the arc on $\partial D_l$ from $t$ to $t_l$. If $K_2\neq \emptyset$ let, for each $j\in K_2$, $w_j\in\partial D_m$ and $x_j\in\partial D_j$ be such that $\overline{w_jx_j}$ is the lower cross tangent segment common to $\partial D_m$ and $\partial D_j$, and let $\delta:= \min \{\delta_j:j\in K_2\}$ and $n:= \min \{j\in K_2:\delta_j=\delta\}$. Let $v'\in \partial D_n$ and $t'\in \partial D_l$ be such that $\overline{v't'}$ is the cross tangent segment common to $\partial D_n$ and $\partial D_l$. If $K_3:=\{j\in K_2:j<n\;\text{and}\;\overline{v't'}\cap D_j\neq \emptyset\}=\emptyset$, take $\Gamma_2:=C\cup \overline{u_mv_m}\cup D'\cup \overline{w_nx_n}\cup E'\cup \overline{v't'}\cup F$, where $D'$ is the arc on $\partial D_m$ from $v_m$ to $w_n$, $E'$ is the arc on $\partial D_n$ from $x_n$ to $v'$ and $F$ is the arc on $\partial D_l$ from $t'$ to $t_l$. If $K_3\neq \emptyset$ we repeat the process which must end. If the curve $\Gamma_2=\Gamma_2^l$ so obtained (which is below $\overline{s_lt_l}$) is inside $G$, replace $\overline{s_lt_l}$ in $\Gamma_1$ by $\Gamma_2^l$; otherwise, some (or all) line segment meet $\partial G$ and, therefore, meet some $D_j$ with $j\in J_1$. This line segments must be corrected with curves of ``third level", and so on. The process allow us to  construct ``polygonal" curves where the slopes are alternative increasing and decreasing. Since the number of $D_i'$,s is finite, the process has an end. Therefore we end up with an admissible arc $\Gamma \subset G$ from $a$ to $b$.\\
	\underline{$J_1=\emptyset $}: In this case either, there is $i\in J$, with $\operatorname{Re}c_i<\operatorname{Re}a$ (if $\arg (b-a)>0$)  such that $\overline{ab}\cap D_i\neq \emptyset$ or, there is $i\in J$, with $\operatorname{Re}c_i>\operatorname{Re}b$ (if $\arg (b-a)<0$) such that $\overline{ab}\cap D_i\neq \emptyset$. In these two cases we first construct $\Gamma_1$ with increasing slope; if $\Gamma_1\subset G$, then take $\Gamma :=\Gamma_1$; otherwise, we correct the line segments with ``polygonal" curves (of second level) with decreasing slope, and so on. Since the process must finish we finally obtain the admissible $\Gamma$ from $a$ to $b$.\\
\textbf{B}: $\overline{ab}$ is vertical. Let $J_1:=\{i\in J: \operatorname{Re}c_i>\operatorname{Re}a\;\text{and}\;D_i\cap \overline{ab}\neq \emptyset\}$.\\ \underline{$J_1\neq \emptyset$}: We first construct the general ``polygonal" $\Gamma_1$ with decreasing slope (to the left of $\overline{ab}$). The rest of the argument is analogous to the previous cases.\\
\underline{$J_1 = \emptyset$}: Here, $K_1:=\{i\in J:\operatorname{Re}c_i<\operatorname{Re}a\;\text{and}\;D_i\cap \overline{ab}\neq \emptyset\}\neq \emptyset$. We first construct the general polygonal $\Gamma_1$ with increasing slope (to the right of $\overline{ab}$) which, if necessary, must be corrected with ``polygonal" curves of ``second level" (with decreasing slope), etc.
\end{proof}
\subsection{Lower bound for the hyperbolic density }\label{sec2:subsec3}
Let $B$ be the horo-crescent domain given by $B=\mathbb{D}\setminus\{z:|z-1/2|\leq 1/2\}$.
We wish to express the density $\nu_B(z)$ in terms of the hyperbolic distance  $d_B(z)$, relative to $\mathbb{D}$, from $z$ to the horocycle $\{z: \vert z-1/2 \vert =1/2\}$. To do this we turn to the conformal invariant model of the hyperbolic plane, namely, the upper half-plane $\mathbb{H}=\{w: \operatorname{Im}w>0\}$ endowed with the metric $\lambda_{\mathbb{H}}(w)|dw|=\dfrac{|dw|}{2\operatorname{Im}w}$, where the computations are sometimes simpler. The Möbius transformation $\tau (z)=i(1+z)/(1-z)$ maps $\mathbb{D}$ conformally onto $\mathbb{H}$ and sends $B$ onto the horizontal strip $=\{w: 0<\operatorname{Im}w < 1\}$. Invariance of the hyperbolic metric under conformal mappings implies that $\nu_B(z)=\lambda_S(\tau (z))/\lambda_{\mathbb{H}}(\tau(z))$. Also the hyperbolic distance is invariant under conformal mappings, therefore $d_B(z)=\delta_S (\tau(z))$, where $\delta_S (\tau(z))$ is the hyperbolic distance, relative to $\mathbb{H}$, from $\tau(z)$ to the line $\operatorname{Im}(w)=1$.
Example (iii) in \cite{Mi}, p. 62, can be used to show that the hyperbolic density of a horizontal strip $\mathcal{S}$, symmetric about the line $\operatorname{Im}(w)=ic, c\in\mathbb{R}$, and width $K\pi$, is given by
\begin{equation}\label{eq-2.3-1}
	\lambda_{\mathcal{S}}(w)=\frac{1}{2K\cos\frac{\operatorname{Im}(w)-c}{K}}.
\end{equation}
(We warn the reader that we have normalized the hyperbolic metric to have Gaussian curvature $-4$). Hence,
\begin{equation}\label{eq-2.3-2}
	\lambda_S(\tau(z))=\frac{\pi}{2}\dfrac{1}{\cos \big [\pi (\operatorname{Im}\tau(z)-1/2)\big ]}=\frac{\pi}{2}\frac{1}{\sin \big (\frac{\pi}{2}\operatorname{Im}\tau(z)\big )}.   
\end{equation}
Now,
\[
\delta_S(\tau(z))=\frac{1}{2}\int_{v}^{1}\frac{dt}{t}=\frac{1}{2}\log\frac{1}{v},
\]
where $v=\operatorname{Im}\tau(z)$. Hence, $\operatorname{Im}\tau(z)=\dfrac{1}{e^{2\delta_S(\tau(z))}}$. This, together with \eqref{eq-2.3-2} yield
\begin{equation}\label{eq-2.3-3}
	\frac{\lambda_S(\tau(z))}{\lambda_{\mathbb{H}}(\tau(z))}=\frac{\pi}{e^{2\delta_S(\tau(z))}}\frac{1}{\sin \big (\frac{\pi}{e^{2\delta_S(\tau(z))}}\big )}.
\end{equation}
Since the hyperbolic distance is conformal invariant, then $d_B(z)=\delta_S(\tau(z))$. So, by \eqref{eq-2.3-3} and the invariance of the hyperbolic density we arrive at the formula
\begin{equation}\label{eq-2.3-4}
	\nu_B(z)=\frac{\pi}{e^{2\delta_B(z)}}\frac{1}{\sin \Big (\frac{\pi}{e^{2\delta_B(z)}}\Big )}.
\end{equation}
The above formula is valid for any horo-crescent domain $\Omega\subset \mathbb{D}$ because any two horo-crescent domains are \textrm{M\"ob}$(\mathbb{D})$-equivalent.
\begin{proof}[\textbf{Proof of Theorem \ref{lower_bound}}]
	Fix $a\in G$. Choose $c\in \partial G$ such that $d_G(a)=d_{\mathbb{D}}(a,c)$. By definition of horo-convexity there exists a horodisk $H$ such that $c \in \partial H$ and $G\cap H=\emptyset$. Let $\Omega$ be the horo-crescent domain $\mathbb{D}\setminus\overline{H}$. Then $G\subset \Omega$. The monotonicity property of the hyperbolic metric yields $\nu_G(a)\geq \nu_{\Omega}(a)$ with equality if and only if $G=\Omega$. Since $d_G(a)=d_{\Omega}(a)$, this inequality in conjunction with \eqref{eq-2.3-4} above completes de proof.	
\end{proof}
The proof of Corollary \ref{corollary_lowerbound} is very similar to the proof given in \cite{Mi}, p. 65 of the corollary to Theorem 4. We give the details just for completeness.
\begin{proof}[\textbf{Proof of Corollary \ref{corollary_lowerbound}}]
	The principle of hyperbolic metric gives $\nu_G(f(z))|f'(z)|/(1-|f(z)|^2)\leq \lambda_{\mathbb{D}}(z)$ for $z\in \mathbb{D}$ with equality if and only if $f$ is a conformal mapping of $\mathbb{D}$ ond $G$. Theorem \ref{lower_bound} then implies that $1/h(d_G(f(z)))\leq \nu_G(f(z)$) with equality if and only if $G$ is $\textrm{M\"ob}(\mathbb{D})$-equivalent to a horo-crescent domain. By combining the two preceding inequalities and the necessary and sufficient conditions for equality, we obtain the corollary.
\end{proof}		

\end{document}